\documentclass[11pt]{article}
\usepackage{amsmath}
\numberwithin{equation}{section}
\usepackage{amsfonts,amssymb,latexsym,amsbsy, bbm, theorem,enumerate,color}
\usepackage{amssymb}
\usepackage{graphicx, graphics}
\usepackage{float,color,fancybox,shapepar,setspace,hyperref}
\usepackage{subfigure}
\usepackage{pgf,tikz}
\usetikzlibrary{arrows}
\makeatletter

\newcommand{\Rmnum}[1]{\expandafter\@slowromancap\romannumeral #1@}
\makeatother
\newtheorem{Main Theorem}{Main Theorem}

\newtheorem{Conjecture}{Conjecture}
\newtheorem{Theorem}{Theorem}

\newtheorem{Lemma}{Lemma}
\newtheorem{Claim}{Claim}

\newtheorem{definition}{Definition}[section]
\def\square{\hbox{\vrule height8pt depth0pt
\vbox{\hrule width7.2pt\vskip7.2pt\hrule width7.2pt}\vrule
height8pt depth0pt}\smallskip}

\def\pf{\medskip\noindent {\emph{\bf Proof}.}~~}
\def\ol{\overline}

\begin{document}

\title{ Tur\'{a}n number  for odd-ballooning of trees}

\author{Xiutao ZHU and Yaojun CHEN\footnote{Corresponding author. Email: yaojunc@nju.edu.cn} \\
 \small{Department of Mathematics, Nanjing University, Nanjing 210093, P.R. CHINA}}
\date{}
\maketitle
\begin{abstract}
The Tur\'an number $ex(n,H)$  is the maximum number of edges in an $H$-free graph on $n$ vertices. Let $T$ be  any tree. The odd-ballooning of $T$, denoted by $T_o$, is a graph obtained by replacing each edge of $T$ with an odd cycle containing the edge, and all new vertices of the odd cycles are distinct. In this paper, we determine the exact value of $ex(n,T_o)$ for sufficiently large $n$ and $T_o$ being good, which generalizes all the known results on $ex(n,T_o)$ for $T$ being a star,
due to Erd\H{o}s et al. (1995), Hou et al. (2018) and Yuan (2018), and provides some counterexamples with chromatic number 3 to a conjecture of Keevash and Sudakov (2004), on the maximum number of edges not in any monochromatic copy of $H$ in a $2$-edge-coloring of a complete graph of order $n$.
\vskip 2mm


\end{abstract}

\section{Introduction}
Let $G=G(V(G), E(G))$ be a graph and $e(G)=|E(G)|$.  For $v\in V(G)$ and $S\subseteq V(G)$, define $N_S(v)=\{u~|~uv\in E(G)~and~u\in S\}$ and $|N_S(v)|=d_S(v)$, and set $N_S(v)=N(v)$ and $|N(v)|=d(v)$ if $S=V(G)$. The {\it maximum degree} of $G$ is  $\Delta(G)=max\{d(v)~|~v\in V(G)\}$.
For $X\subseteq V(G)$,  $G[X]$ denotes the subgraph of $G$ induced by $X$.  If $X,Y\subseteq V(G)$ and $X\cap Y=\emptyset$, then $G[X,Y]$ denotes the induced bipartite  subgraph of $G$ with bipartition $(X,Y)$. Let $G,H$ be two graph, we use $G\cup H$ and $G+H$ to denote the disjoint union and join of $G$ and $H$, respectively. A path, cycle, complete graph, empty graph and star on $n$ vertices are denoted by $P_n$, $C_n$, $K_n$, $E_n$ and $S_{n-1}$, respectively. For a family of graphs $\mathcal{H}$, $G$ is called $\mathcal{H}$-{\it free} if $G$ contains no member of $\mathcal{H}$ as a subgraph.  Write $H$-free for $\mathcal{H}$-free if $\mathcal{H}=\{H\}$. A {\it covering} of $G$ is a set of vertices which  meets all edges. Let $\beta(G)$ denote the minimum number of vertices in a covering of $G$. Define  $\nu(G)$  as the size of a maximum matching of $G$. Set
$$f(\nu,\Delta)=max\{e(G):\nu(G)\le \nu~\text{and}~\Delta(G)\le \Delta\}.$$
The following is a celebrated theorem due to Chv\'{a}tal and Hanson.
\begin{Theorem}(Chv\'{a}tal and Hanson \cite{Chvatal})\label{sunflower}
$$ f(\nu, \Delta)=\nu\Delta +\left\lfloor\frac{\Delta}{2}\right\rfloor\left\lfloor\frac{\nu}{\lceil\Delta/2\rceil}\right\rfloor\le \nu(\Delta+1).$$
\end{Theorem}
The special case when $\nu=\Delta=k-1$ was first determined by Abbott et al. \cite{AHS}:

$$f(k-1,k-1)=\begin{cases}
       k^2-k &if~k~is~odd,\\
        k^2-\frac{3}{2}k &if~k~is~even.
    \end{cases}$$

The Tur\'an number of $H$, denoted by $ex(n,H)$, is the maximum size in an $H$-free graph on $n$ vertices. Extremal graph theory dates back to the early 1940's when Tur\'an \cite{turan} proved that $T_p(n)$ is the unique extremal graph of $ex(n,K_{p+1})$, where $T_p(n)$  is a complete $p$-partite graph on $n$ vertices in which each partite set has $\lfloor n/p\rfloor$ or $\lceil  n/p\rceil$ vertices, called {\it Tur\'an graph}.
The famous Erd\H{o}s-Stone-Simonovits Theorem \cite{Er,stone} states if $H$ is a graph with chromatic number $\chi(H)=p+1\ge 3$, then
$$ex(n,H)=e\big(T_p(n)\big)+o(n^2).$$
This means  the Tur\'an number $ex(n,H)$ is determined asymptotically for $H$ being a nonbipartite graph. Moreover, Erd\H{o}s and Simonovits \cite{Er} showed that   the ``asymptotic structure" of the extremal graphs for the Tur\'an number $ex(n,H)$ is also determined by the chromatic number $\chi(H)$ of $H$.  This is to say, the extremal graphs for a nonbipartite $H$ is very similar to $T_p(n)$ when $n$ is sufficiently large. For convenience to describe the ``asymptotic structure" of the extremal graphs for $ex(n,H)$,  Simonovits \cite{Simonovits} defined the following decomposition family.
\begin{definition}
The decomposition family $\mathcal{M}(H)$ is the set of minimal graphs $M$ such that if  an $M$ is embedded into one partite set of $T_p(n)$, then
the resulting graph contains $H$ as a subgraph, where $\chi(H)=p+1$ and $n$ is large enough.
\end{definition}
This concept provides us an idea on how to consider the extremal graphs for $ex(n,H)$ via Tur\'an graph $T_p(n)$, that is, one can embed a maximal $\mathcal{M}(H)$-free graph into one partite set of $T_p(n)$ to obtain the extremal graphs for $ex(n,H)$.
Although  it is still a challenging problem to  determine the  Tur\'an number and all extremal graphs for many nonbipartite graphs, Simonovits' method is still a very useful tool for dealing Tur\'an problems.
The following result, on the Tur\'an number of friendship graph consisting of $k$ triangles intersecting in one common vertex, due to Erd\H{o}s et al., can be viewed as a classical example using this idea.
\begin{Theorem}(Erd\H{o}s et al. \cite{friendship})\label{friendship}
For any $k\ge 1$ and $n\ge 50k^2$,
$$ex(n,K_1+kK_2)=ex(n,K_{3})+f(k-1,k-1).$$
\end{Theorem}

 Note that the decomposition family of friendship graph is $\{kK_2,S_k\}$ and the maximum number of edges of a $\{kK_2,S_k\}$-free graph is exactly $f(k-1,k-1)$. It was also shown in \cite{friendship} that the only extremal graph in Theorem \ref{friendship}  is the one obtained from $T_2(n)$ by embedding a maximal $\{kK_2,S_k\}$-free graph into one partite set.

Theorem \ref{friendship} was generalized in different ways. Since a triangle is a $K_3$, Chen et al. \cite{Chen}  considered the same problem by replacing each triangle with a $K_{p+1}$, $p\geq 2$, in Theorem \ref{friendship}, and obtained a generalization as follows.

\begin{Theorem}(Chen et al. \cite{Chen})
For any $p\ge 2$ and $k\ge 1$, when $n\ge 16k^3(p+1)^8$,
$$ex(n,K_1+kK_{p})=ex(n,K_{p+1})+f(k-1,k-1).$$
\end{Theorem}
 Inspired by these results, Glebov \cite{Glebov} considered the same extremal problem for many  $K_{p+1}$ intersecting in a different way: replace each edge of $P_{k}$ with a $K_{p+1}$.
Short after that, Liu \cite{H.Liu} introduced the concept of edge blow-up graph. Let $H$ be a graph and $p$ an integer, the {\it edge blow-up} of $H$, denoted by $H^{p+1}$, is one obtained by replacing each edge in $H$ with a clique $K_{p+1}$, where the new vertices are all different. Clearly, the graphs mentioned above are all special edge blow-up graphs: $S_k^3$, $S_k^{p+1}$ and $P_{k}^{p+1}$. Liu \cite{H.Liu} also determined the Tur\'an number for the edge blow-up of some trees and all cycles.
Recently, Wang et al. \cite {Wang} determined the extremal graphs for the edge blow-up of a large family of trees and Yuan \cite{Yuan} got a tight bound for any $G^{p+1}$ when $p\ge\chi(G)$.

Observe a triangle is also a $C_3$, Hou et al. \cite{Hou,Hou2} and Yuan \cite{Yuan2} started to consider generalizing Theorem \ref{friendship} in another direction: Replace each $C_3$ with an odd cycle. Let $H$ be any graph. The {\it odd-ballooning} of $H$  is a graph obtained by replacing each edge of $H$ with an odd cycle containing the edge, and all new vertices of the odd cycles are distinct. Obviously, the friendship graph is a special odd-ballooning of a star.
Hou et al. \cite{Hou2} and Yuan \cite{Yuan2} studied the Tur\'an number of any odd-ballooning of a star and got the following.
\begin{Theorem}(Hou et al. \cite{Hou2} and Yuan \cite{Yuan2})
Let $F_{k,k_1}$ be the graph consisting of $k$ odd cycles intersecting in one common vertex and $k_1$ of these cycles are triangles. When $n$ is large enough,
$$ex(n,F_{k,k_1})= \left\lfloor\frac{n^2}{4}\right\rfloor+\begin{cases}
     (k-1)^2 & {\text if}~k> k_1, \\
       f(k-1,k-1)& {\text if}~k=k_1.
    \end{cases}$$
\end{Theorem}
Later, Zhu et al. \cite{Zhu} determined the Tur\'an number for the odd-ballooning of a path.

No matter a star or a path, it is a special tree. This motivates us to consider the Tur\'an number for the odd-ballooning of a tree $T$ in a more general situation.
Before stating our result, we first introduce some additional notations.
Throughout this paper, we use $T=T[A,B]$ to denote a tree $T$ with unique bipartition $A$, $B$, and assume $a=|A|\le |B|$ and $\delta(A)=min\{d(v)~|~v\in A\}$.  An edge $uv$ is called a {\it leaf-edge} if at least one end of $uv$is of degree $1$.  The odd-ballooning of $T$ is denoted by $T_o$. For any $uv\in E(T)$,  we use $C(u,v)$ to denote the odd cycle containing $uv$ in $T_o$ and $P(u,v)=C(u,v)-uv$. If $C(u,v)$ is a triangle in $T_o$,  then we say the edge $uv$ is of Type \Rmnum{1} and $uv$ is of Type \Rmnum{2} otherwise. An odd-ballooning $T_o$ of $T$ is {\it good}  if all  edges of Type I are {\it leaf-edges} and the non-leaf vertices  are in $A$. For $u\in A$, we use $d_I(u)$ to denote the number of  the edges of Type \Rmnum{1}  incident with $u$.

We still use Simonovits' idea to investigate the Tur\'an number $ex(n,T_o)$ and the extremal graphs.
Let $G(n,p,a)=E_{a-1}+T_p(n-a+1)$,  $X=V(E_{a-1})$ and $X_1,...,X_p$ be the $p$ partite sets of $T_p(n-a+1)$. We use $G(n,p,a,F)$ and $G(n,p,a,\mathcal{H},F)$ to denote  a graph  obtained from $G(n,p,a)$ by embedding a graph $F$ into $X_1$, and embedding a maximum $\mathcal{H}$-free graph into $X$ and  a graph $F$ into $X_1$, respectively.
Moreover, we define a family of subgraphs $\mathcal{B}(T_o)$ based on the decomposition family $\mathcal{M}(T_o)$ as follows.

\begin{definition}\label{B} If each $M\in \mathcal{M}(T_o)$ has no covering with less than $a$ vertices, then $\mathcal{B}(T_o)=\{ K_a \} $ and otherwise,
$$\mathcal{B}(T_o)=\{ M[S]~|~M\in \mathcal{M}(T_o)~and~S~is~a~covering~of~M~with~|S|<a  \}.$$
\end{definition}

The main result of this paper is the following.
\begin{Theorem}\label{Thm1}
Let $T=T[A,B]$  be a tree with $a=|A|\leq |B|$ and $T_o$ a good odd-ballooning of $T$. Suppose $u\in A$ with $d(u)=\delta(A)=k$ and $d_I(u)=k_1$. If $n$ is large enough, then
\[
ex(n, T_o)=e(G(n,2,a))+ex(a-1,\mathcal{B}(T_o))+\begin{cases}
(k-1)^2, & if ~k> k_1,\\
f(k-1,k-1), & if~k=k_1.
\end{cases}
\]
Furthermore, if $k_1< k$, then $G(n,2,a,\mathcal{B}(T_o), K_{k-1,k-1})$ is an extremal graph, and if $k_1=k$, then $G(n,2,1, F)$ is an extremal graph, where $F$ is the extremal graph of the function $f(k-1,k-1)$.

\end{Theorem}

Since a star $S_\ell$ is a tree $T[A,B]$ with $|A|=a=1$ and any odd-ballooning of $S_\ell$ is good,  one can see
Theorem \ref{Thm1} generalizes the results obtained in \cite{friendship,Hou,Hou2,Yuan2}.
Moreover, let $f(n, H)$ denote the maximum number of edges not in any monochromatic copy of $H$ in a $2$-edge-coloring of $K_n$. As a by-product, Theorem \ref{Thm1} also provides some counterexamples with chromatic number 3 to the following conjecture:
 \begin{Conjecture}(Keevash and Sudakov \cite{Sudakov})\label{conj1}
For any $H$ and sufficiently large $n$, $$f(n,H)= ex(n,H).$$
\end{Conjecture}
The details will be presented in Section 5. The other parts of this paper is organized as follows. Section 2 is devoted to characterizing the decomposition family of $T_o$;   Section 3 contains some preliminaries; The proof of Theorem \ref{Thm1} is given in Section 4.

\section{Decomposition family of $T_o$}
In this section, our task is to characterize the decomposition family of $T_o$ through two operations. Let $G$ be a graph, $u\in V(G)$ and $N(u)=\{v_1,\ldots,v_d\}$. A {\it splitting} on $u$ is to replace $u$ with an independent set $\{u_1,\ldots,u_d\}$ and  $uv_i$ with a new edge $u_iv_i$, $1\leq i\leq d$, and write  $s^{-1}(u_iv_i)=uv_i$. Let $uv$ be a leaf-edge. If $d(u)=1$ and $d(v)\geq 2$, then {\it peeling off} $uv$ is to delete $uv$ and add a new vertex $v'$ and a new edge $uv'$, and write $p(uv)=uv'$. If $d(u)=d(v)=1$, then {\it peeling off} $uv$  means do nothing and $p(uv)=uv$.

For a given odd-ballooning $T_o$ of a tree $T$, we use $\mathcal{SP}(T)$ to denote the family of graphs, each of which can be obtained from $T$ through splitting the vertices in some independent set first, say the resulting graph $T_1$. Then peeling off some  leaf-edges of $T_1$ which satisfies  $uv$ or $s^{-1}(uv)$ is of Type II.

We have the following lemma.

\begin{Lemma}\label{lemma1}
 For any tree $T$ and any odd-ballooning $T_o$,  $\mathcal{M}(T_o)=\mathcal{SP}(T)$. In particular, if $T_o$ is good, then a matching of size $e(T)$ is in $\mathcal{M}(T_o)$.
\end{Lemma}
\pf
 Let $X,Y$ be two partite sets of $T_2(n)$, where $n$ is sufficiently large.

 Let $M$ be any graph in $\mathcal{M}(T_o)$. We first show $M\in\mathcal{SP}(T)$.
 Embed $M$ into $X$, then  by the definition of $\mathcal{M}(T_o)$, the resulting graph contains a copy of $T_o$. Color the vertices in $X$ red and  the vertices in $Y$ blue, and call an edge red if its two ends are colored red. Clearly, $E(M)$ are all red edges.

Because $\chi(T_o)=3$ and $\chi(T_o-E(M))=2$, hence $|E(M)\cap E(C(u,v))|\ge1$. Suppose $C(u,v)$ is an odd cycle in $T_o$ containing at least two red edges. Since $X$, $Y$ are large enough, if $uv$ is a red edge, then  we can replace $P(u,v)=ua_1\cdots a_{i}v$ by a new proper colored path $ub_1\cdots b_{i}v$ using vertices distinct with the original $T_o$, and if $u$ is red and $v$ is blue, then we can replace $P(u,v)$ by a path which contains exactly one red edge of $P(u,v)$. Obviously, the red edges in the new $T_o$ form a subgraph of $M$, contradicting the minimality of $M$.  Therefore, $C(u,v)$ contains exactly one red edge.

Now, consider the skeleton $T$ of $T_o$. The blue vertices  in $T$ (if any) is an independent set. Split all blue vertices of $T$.
Suppose $u$ is any blue vertex, $N_T(u)=\{v_1,\ldots, v_d\}$ and $u$ is split into $\{u_1,\ldots,u_d\}$. Clearly, $u_iv_i$  is a {\it leaf-edge} after splitting. Let $ww'$ be the only red edge in $C(u,v_i)$. If $ww'$ is not incident with $v_i$, then $C(u,v_i)$ is of order at least $5$ and $uv_i$ is of Type \Rmnum{2}. In this case, we peel off $u_iv_i$ and let $p(u_iv_i)$ be the new edge. If $ww'$ is incident with $v_i$, then we do not peel off $u_iv_i$. Let $T'$ be the resulting graph by peeling off all such edges. Then $T'\in \mathcal{SP}(T)$. Now, let $e\mapsto  e$ if $e\in E(M)\cap E(T')$, $ww'\mapsto p(u_iv_i)$ if $ww'$ is not incident with $v_i$ and $ww'\mapsto u_iv_i$ if $ww'$ is incident with $v_i$, we can see that $M\cong T'$, and so $M\in\mathcal{SP}(T)$.

On the other hand, let $T'$ be any graph in $\mathcal{SP}(T)$, which is obtained from $T$ by splitting some independent set $U$ first, and then peeling off some {\it leaf-edges} satisfying $uv$ or $s^{-1}(uv)$ is of Type \Rmnum{2}, from the resulting graph. We will show $T' \in\mathcal{M}(T_o)$.

For any $u\in V(T)$, color $u$ with blue if $u\in U$ or $uv\in E(T)$ is a {\it leaf-edge} with $d(u)=1$ and $v\not\in U$, otherwise color $u$ with red. And then, color the vertices in $V(T_o)- V(T)$ as follows. Let $uv$ be any edge of $T$. If $uv$ is red, then give a proper red-blue coloring to $P(u,v)$; If $u$ is blue and $v$ is red,
then $d_T(u)=1$ and $v\notin U$, or $u\in U$ is split and say $s(uv)=u'v$. In this case, give $P(u,v)$ a red-blue coloring such that it contains exactly one red edge $ww'$, and $ww'$ is incident with $v$ if $uv$ or $u'v$ is not peeled off and not incident $v$ if $uv$ or $u'v$ is peeled off.
Because if $uv$ or $u'v$ is peeled off, then $uv$ is of Type \Rmnum{2} and $P(u,v)$ is an odd path of order at least $5$, and so such a coloring exists and all blue vertices form an independent set.
Assume that $M'$ is the subgraph in $T_o$ induced by all red edges.
Let $e\mapsto  e$ if $e\in E(T)\cap E(M')$, $p(uv)\mapsto ww'$ if $d_T(u)=1$ and $uv$ is peeled off, $u'v\mapsto ww'$ if $ww'$ is incident with $v$ and $p(u'v)\mapsto ww'$ if $ww'$ is not incident with $v$, we can see that $T'\cong M'$.

 Observe that each odd cycle in $T_o$ contains exactly one red edge, we can see that $T_o-E(M')$ is proper red-blue colored. Thus, we can embed $T_o-E(M')$ into $T_2(n)$ such that all red vertices of $T_o$ are in $X$ and all blue vertices of $T_o$ are in $Y$. Note that $T'\cong M'$ and $e(T')=e(T)$, we have $T'\in \mathcal{M}(T_o)$.

\vskip 2mm
In particular, if $T_o$ is a good odd-ballooning of $T=T[A,B]$, then apply vertex splitting on $A$ first, the resulting graph is the disjoint union of stars and all edges of Type \Rmnum{1} become isolated edges. For any star other than a $K_2$ in the resulting graph, each edge $uv$ of it satisfying $uv$ or $s^{-1}(uv)$ is
of  Type \Rmnum{2}. Thus, peel off some edges from each such star,  we can obtain a matching of size $e(T)$. Since $\mathcal{M}(T_o)=\mathcal{SP}(T)$, this matching is in
$\mathcal{M}(T_o)$.$\hfill\blacksquare$

\section{Preliminaries}

\begin{Lemma}(K\"{o}nig \cite{Konig})\label{konig}
Let $G$ be a bipartite graph, then $\beta(G)=\nu (G)$.
\end{Lemma}

\begin{Lemma}(Hall \cite{Hall})\label{Hall} Let $G=G[X,Y]$ be a bipartite graph. Then $\nu(G)\ge |X|$ if and only if
$$|N(S)|\ge |S|\quad \text{for all} \quad S\subseteq X.$$
\end{Lemma}

\begin{Lemma}(Wang et al. \cite{Wang})\label{Hou}
Let $T[A,B]$ be a tree with $\delta(A)=k\ge 2$. If a vertex in $A$ is split, then the resulting graph $T'$ satisfies $\nu(T')\ge a-1+k$.
\end{Lemma}

\begin{Lemma}\label{lemma2}
Let $T=T[A,B]$ be a tree and $T_o$ be any odd-ballooning of  $T$. Then $\mathcal{B}(T_o)=\{K_{a}\}$ if and only if $\beta(T)=a$. Furthermore,  if $\delta(A)\ge 2$, then $\beta(T)=a$.
\end{Lemma}
\pf By Lemma \ref{lemma1}, $\mathcal{M}(T_o)=\mathcal{SP}(T)$. Let $T'$ be any graph in $\mathcal{SP}(T)$. Since splitting vertices and peeling off {\it leaf-edges}
 do not decrease the size of maximum matching, we have $\nu(T')\geq \nu(T)$.
 If $\beta(T)=a$,  then by Lemma \ref{konig}, $\beta(T')=\nu(T')\ge \nu(T)=\beta(T)=a$.
 That is, $T'$ has no  covering $S$  with $|S|<a$. By Definition \ref{B}, we have $\mathcal{B}(T_o)=\{K_{a}\}$. On the other hand, if $\mathcal{B}(T_o)=\{K_{a}\}$, then since $T\in \mathcal{M}(T_o)$, we get $\beta(T)=a$.

Furthermore, if $\delta(A)\ge 2$, then since $T[S,N(S)]$ is a forest for any $S\subseteq A$, we have
$$2|S|\le e(T[S,N(S)])\le |S|+|N(S)|-1,$$
which implies $|N(S)|\ge |S|+1$. By Lemma \ref{Hall},  we have $\beta(T)=\nu(T)=a$ $\hfill\blacksquare$

\begin{Lemma}\label{star-matching} Let $G$ be an $\{S_k,kK_2,S_{k-1}\cup K_2\}$-free graph without isolated vertices. Then $e(G)\leq (k-1)^2$ with equality if and only if $G=K_{k-1,k-1}$, or $G=3K_3$ and $k=4$.
\end{Lemma}
\pf Let $G$ be an $\{S_k,kK_2,S_{k-1}\cup K_2\}$-free graph with $e(G)\ge (k-1)^2$. Obviously,  $\Delta(G)\le k-1$. If $\Delta(G)\le k-2$, then $G$ is $\{kK_2,S_{k-1}\}$-free and hence we have
$$(k-1)^2\le e(G)\le f(k-1,k-2)=(k-1)(k-2)+\left\lfloor\frac{k-2}{2}\right\rfloor\left\lfloor\frac{k-1}{\lceil(k-2)/2\rceil}\right\rfloor,$$
which implies $k=4$ and $G=3K_3$.
If $\Delta(G)=k-1$, let $v$ be a vertex with $d(v)=k-1$. Since $G$ is $S_{k-1}\!\cup\! K_2$-free, $G-N(v)$ has no edges, and so $e(G)\le (k-1)^2$, with equality  if and only if all vertices in $N(v)$ have degree $k-1$ and $N(v)$ is independent set. Moreover, all vertices in $N(v)$ have the same neighborhoods for otherwise we can find an $S_{k-1}\cup K_2$ in $G$. Therefore, $G=K_{k-1,k-1}$.$\hfill\blacksquare$

\begin{Lemma}(Erd\H{o}s et al. \cite{friendship})\label{lemma3}
Let $\Delta$ and $b$ two nonnegative integers such that $b\le \Delta-2$. If $\Delta(G)\le \Delta$ , then
$$\sum_{v\in V(G)}min\{d(v),b\}\le \nu(G)\big(b+\Delta\big).$$
\end{Lemma}

\vskip 2mm
Suppose $G$ is a graph with partition $V(G)=V_0\cup V_1$, let $G_0=G[V_0]$, $G_1=G[V_1]$, $G_{cr}=G[V_0,V_1]$, and $d_{cr}(v)$ be the degree of the vertex $v$ in $G_{cr}$. We have the following.

\begin{Lemma}\label{lemma4}
Let $k\ge k_1$ be two nonnegative integers. If the following hold,
\vskip 2mm
(1) $G_i$ is $\{S_{k-\ell}\cup \ell K_2 : 0\leq \ell\leq min\{k-k_1,k-2\}~or~\ell=k-1\}$-free for $i=0,1$,


(2) $d_{V_i}(v)+\nu(G_{1-i}[N_{V_{1-i}}(v)])\le k-1$ for any $v\in V_i$,

(3) If $k>k_1$,  then $N_{V_{1-i}}(v)$ are isolated  in $G_{1-i}$ for any $v\in V_i$ with $d_{V_i}(v)=k-1$,
\vskip 2mm
\noindent then
\begin{equation}\label{eq3.1}
e(G_0)+e(G_1)-\big(|V_0||V_1|-e(G_{cr})\big)\le\begin{cases} (k-1)^2 &{\text for~} k>k_1,\\
f(k-1,k-1) &{\text for~} k=k_1.
\end{cases}
\end{equation}
\end{Lemma}

\pf  As it is shown in \cite{friendship} that (\ref{eq3.1})  holds for $k=k_1$, we may assume that $k>k_1$.

Choose a graph $G$  with partition $V(G)=V_0\cup V_1$ satisfying (1), (2) and (3), such that  $e(G_0)+e(G_1)-\big(|V_0||V_1|-e(G_{cr})\big)$ is as large as possible, and subject to this, $|G|$ is as small as possible.

If there is some $v\in V_0$ such that $d_{V_0}(v)-\big(|V_1|-d_{cr}(v)\big)\leq 0$, then since $G-v$ with partition $V(G-v)=(V_0-\{v\})\cup V_1$ still satisfies (1), (2) and (3),  $|G-v|<|G|$ and
\[\begin{split}
& ~e(G_0-v)+e(G_1)-\big(|V_0-\{v\}||V_1|-e(G_{cr}-v)\big)\\
\geq &~ e(G_0)+e(G_1)-\big(|V_0||V_1|-e(G_{cr})\big),
\end{split}\]
this contradicts the choice of $G$. Thus,
we have $d_{V_i}(v)-\big(|V_{1-i}|-d_{cr}(v)\big)>0$ for any $v\in V_i$ by the symmetry of $V_0$ and $V_1$.
Moreover, because $d_{V_0}(v)+\nu(G_1[N_{V_1}(v)]\le k-1$, we have   $d_{V_0}(v)-\big(|V_1|-d_{cr}(v)\big)
\le k-1-\big(\nu(G_1[N_{V_1}(v)])+|V_{1}|-d_{cr}(v)\big)
 \le  k-1-\nu(G_{1})$.
 The second inequality holds since any matching in $G_1$ has at most $\nu(G_1[N_{V_1}(x)])$ edges in $N_{V_1}(v)$ and at most $|V_1|-d_{cr}(v)$ additional edges.
 Therefore, for any $v\in V_i$, $i=0,1$, we have
\begin{equation}\label{3.2}
0<d_{V_i}(v)-\big(|V_{1-i}|-d_{cr}(v)\big)\le k-1-\nu(G_{1-i}).
\end{equation}
By  (\ref{3.2}), we can deduce that
\begin{equation}\label{3.3}
\nu(G_i)\le k-2~for~i=0,1.
\end{equation}

If $\nu(G_i)=0$ for some $i\in \{0,1\}$, then $e(G_i)=0$. Observe that $|V_0||V_1|-e(G_{cr})=e\big(\ol{G}[V_0,V_1]\big)\geq 0$,
 we get $e(G_0)+e(G_1)-\big(|V_0||V_1|-e(G_{cr})\big)\le e(G_{1-i})\le(k-1)^2$ by Lemma \ref{star-matching}, and hence  (\ref{eq3.1}) holds.

If $\nu(G_i)=1$ for some $i\in \{0,1\}$, we assume  that $\nu(G_1)=1$ by symmetry. Clearly, $G_1$ is a star or a triangle with some isolated vertices and $e(G_1)\le max\{k-1,3\}$. Since $G_0$ is $S_k$-free, $\Delta(G_0)\leq k-1$. If $\Delta(G_0)=k-1$, then $e(G_0)\le f(k-2,k-1)\le (k-2)k$ by (\ref{3.3}) and Lemma \ref{sunflower}. Let $v\in V_0$ with $d_{V_0}(v)=k-1$. By the assumption (3), $N_{V_1}(v)$ are isolated vertices of $G_1$, which implies the vertices of the star or the triangle in $G_1$ are nonadjacent to $v$, and hence
$|V_0||V_1|-e(G_{cr})\ge e(G_1)$. Thus, we have
$$e(G_0)+e(G_1)-\big(|V_0||V_1|-e(G_{cr})\big) \le (k-2)k+e(G_1)-e(G_1)< (k-1)^2.$$
If $\Delta(G_0)\le k-2$, then $e(G_0)\le f(k-2,k-2)$ by Lemma \ref{sunflower} and
$$e(G_0)+e(G_1)-\big(|V_0||V_1|-e(G_{cr})\big)\le f(k-2,k-2)+k-1 \le (k-1)^2.$$

Now, assume that $2\le \nu(G_i)\le k-2$ for each $i$. By (\ref{3.2}), we have
\[\begin{split}
2e(G_i)-\big(|V_0||V_1|-e(G_{cr})\big)=&\sum_{v\in V_i}\big\{d_{V_i}(v)-\big(|V_{1-i}|-d_{cr}(v)\big)\big\}
\\ \le&\sum_{v\in V_i}min\big\{d_{V_i}(v), k-1-\nu(G_{1-i})\big\}.
\end{split}\]
Applying Lemma \ref{lemma3} on $G_i$ with $\Delta=k-1$ and $b=k-1-\nu(G_{1-i})\le \Delta-2$, we have

$$2e(G_i)-\big(|V_0||V_1|-e(G_{cr})\big) \le \nu(G_i)\big(2(k-1)-\nu(G_{1-i})\big),$$
and then
\[\begin{split}
&~ 2e(G_0)+2e(G_1)-2\big(|V_0||V_1|-e(G_{cr})\big)
\\ \le &~  \nu(G_0)\big(2(k-1)-\nu(G_1)\big)+\nu(G_1)\big(2(k-1)-\nu(G_0)\big)
\\ =&~ 2\big(k^2-2k+1-(k-1-\nu(G_0))(k-1-\nu(G_1))\big)<  2(k-1)^2,
\end{split}\]
the last inequality follows from (\ref{3.3}). Thus, we finish the proof of Lemma \ref{lemma4}. $\hfill\blacksquare$

\section{The proof of Theorem \ref{Thm1}}
\subsection{Lower bound of $ex(n, T_o)$}
It is not difficult to check the sizes of the three graphs described in Theorem \ref{Thm1} is the expected value of $ex(n,T_o)$. So it suffices to show each of the three graphs is $T_o$-free to get the lower bound for $ex(n,T_o)$.

 If $d_I(u)=k$, then $T=S_k$ and $T_o$ consists of $k$ triangles intersecting in one common vertex. By Lemma \ref{lemma1}, $\mathcal{M}(T_o)=\{S_k,kK_2\}$.  Let $F$ be an extremal graph of the function $f(k-1,k-1)$, so $F$ is $\{S_k,kK_2\}$-free. Note that embedding an $F$ into one class of $T_2(n)$ is $G(n,2,1,F)$, by the definition of decomposition family, it is $T_o$-free. Hence we may assume $d_I(u)<k$.

 It remains to show $G=G(n,2,a,\mathcal{B}(T_o),K_{k-1,k-1})$ is $T_o$-free.
Since $G$ is obtained by embedding a $K_{k-1,k-1}$ into the class $X_1$  and a maximum $\mathcal{B}(T_o)$-free graph into the class $X$ of $G(n,2,a)$,  by the definition of decomposition family, we need only  to prove $G[X\cup X_1]=G[X]+G[X_1]$ is $\mathcal{M}(T_o)$-free.

If $k=1$, then $e(G[X_1])=0$. Since $G[X]$ is $\mathcal{B}(T_o)$-free, by the definition of $\mathcal{B}(T_o)$, we know $G[X\cup X_1]$ is $\mathcal{M}(T_o)$-free, and so we may assume $k\ge 2$. In this case, we have $G[X]=K_{a-1}$ by Lemma \ref{lemma2}.

Suppose to the contrary that $G[X\cup X_1]$ contains a $T' \in \mathcal{M}(T_o)$. Choose $T'$ such that $|T'|$ is minimum.
By Lemma \ref{lemma1}, $T'\in \mathcal{SP}(T)$. Assume that $T'$ is obtained from $T$ by splitting an independent set $U$ and then peeling off some {\it leaf-edges} satisfying  $uv$ or $s^{-1}(uv)$ is of Type II. Clearly, $T'$ is bipartite. Let $A', B'$ be two partite sets of $T'$ such that if $u\in A$ or $u$ is split from a vertex in $A$, then $u\in A'$,  and the same is for $B$ and $B'$. Such a bipartition $A',B'$ is unique.
Let $A'_1=X\cap A'$, $ A'_2=X_1\cap A'$, $B'_1=X\cap B'$, $B'_2=X_1\cap B'$ and $T'[A'_2,B'_2]$ be the subgraph of $T'$ induced by $A'_2\cup B'_2$.

 If some vertex in $A$ is split, then $\nu(T')\ge a-1+k$ by Lemma \ref{Hou}. Since $A'_1\cup B'_1$ is a vertex covering of $T'-E(T'[A'_2,B'_2])$, then $\nu(T'-E(T'[A'_2,B'_2]))\le |A'_1\cup B'_1| \le a-1$ by Lemma \ref{konig}. Note that $\nu(T'[A'_2,B'_2])\le k-1$ since only $ K_{k-1,k-1}$ is in $X_1$, we have
$$ \nu(T') \le \nu(T'-T'[A'_2,B'_2])+\nu(T'[A'_2,B'_2])\le a+k-2,$$
a contradiction. Thus, no vertex in $A$ is split.

Assume that $uv$ is a {\it leaf-edge} with $d(v)=1$, which is peeled off after splitting $U$, and $p(uv)=u'v$.
Because $\delta(A)=k\ge2$ and no vertex in $A$ is split, we have $u\in A$ and $u'\in A'$. Suppose that we get $T'$ by peeling off $uv$ from $T_*'$. Then $T_*'\in \mathcal{SP}(T)$. We now show $T_*'\subseteq G[X\cup X_1]$.
If $u\in A'_1$, or $u\in A'_2$ and one of $u'$, $v$ is in $X$, then $uv$ or  $uu'$ is an edge of $G[X\cup X_1]$. Add $uv$ or $uu'$, and delete $u'v$, we can find a $T_*'$ in  $G[X\cup X_1]$.
If $u,u'\in A'_2$ and $v\in B'_2$, then since $G[X_1]$ is a $K_{k-1,k-1}$ with some isolated vertices,  one of $u'$, $v$ is adjacent to all $N_{T'}(u)$. Replace $u$ with $u'$ or $v$, we can find a $T_*'$ in  $G[X\cup X_1]$. By Lemma \ref{lemma1},  $T_*'\in \mathcal{M}(T_o)$, which contradicts the choice of $T'$ since $|T_*'|<|T'|$. Thus, no edge is peeling off after splitting $U$.

By the argument above, $A=A'=A'_1 \cup A'_2$. Since $a=|A'_1|+|A'_2| \ge |A'_1|+|B'_1|+1$ and $T'[A'_2, B'_1]$ is a forest, we have
$$e(T'[A'_2,B'_2])\ge k|A'_2|-(|A'_2|+|B'_1|-1)\ge (k-2)|A'_2|+2.$$
If $k=2$, then $e(T'[A'_2,B'_2])\ge 2$, which contradicts that $e(G[X_1])=1$, and hence $k\ge 3$.  Because $\Delta(G[X_1])=k-1$, by the inequality above, $A_2'$ has two vertices of degree $k-1$ in $T'[A'_2,B'_2]$, that lie in the same partite set of $K_{k-1,k-1}$. Thus, the two vertices have  the same neighbors in $B'_2$, and so $T'[A'_2,B'_2]$ has cycles, a contradiction.


Therefore,  $G[X]+G[X_1]$ is $\mathcal{M}(T_o)$-free, and so $G$ is $T_o$-free.

\subsection{ Upper bound of $ex(n, T_o)$}
In order to establish the upper bound for $ex(n,T_o)$, we need a result of Simonovits.
\begin{definition}
Denote by $\mathcal{D}(n,p,r)$ the family of $n$-vertex graphs $G$ satisfying the following symmetry condition:
\begin{itemize}
  \item  It is possible to omit at most $r$ vertices of $G$ so that the remaining graph $G'$ is a join of graphs of almost equal order: $G'=G_1+\cdots+G_p$, where $\left||V(G_i)|-\frac{n}{p}\right|\le r$ for all $i\le p$.
  \item For each $i\le p$,  there exist connected graphs $H_i$ such that $G_i= k_iH_i$, where $k_i=\frac{|V(G_i)|}{|V(H_i)|}$,
and any two copies $H_i',H_i''$ of $H_i$ in $G_i$, are symmetric subgraphs of $G$: there exists an isomorphism $\phi: V(H_i')\mapsto V(H_i'')$ such that for any $u\in V(H_i')$ and $v \in V(G)-V(G')$, $uv\in E(G)$ if and only if $\phi(u)v \in E(G)$.
\end{itemize}
\end{definition}

The graphs $H_i$ ($1\leq i\leq p$) will be called the {\it blocks} and the vertices in $V(G)-V(G')$ will be called {\it exceptional vertices}.

\begin{Theorem}(Simonovits \cite{Simonovits})\label{Simonovits}
 For a given graph $H$ with $\chi (H)=p+1$, if $\mathcal{M}(H)$ contains a linear forest, then there exist $r=r(H)$ and $n_0=n_0(r)$ such that $\mathcal{D}(n,p,r)$ contains an extremal graph of $H$ for all $n\ge n_0$. Furthermore, if this is the only extremal graph in $\mathcal{D}(n,p,r)$, then it is the unique extremal graph for every sufficiently large $n$.
\end{Theorem}

We now begin to show the upper bound of $ex(n,T_o)$.
\vskip 2mm
Let $T_o$ be any good odd-ballooning of a tree $T=[A,B]$ with $a=|A|\leq |B|$. By Lemma \ref{lemma1},  $\mathcal{M}(T_o)$ contains a linear forest $e(T)\cdot K_2$.  By Theorem \ref{Simonovits}, $\mathcal{D}(n,2,r)$ contains an extremal graph $G$ for sufficiently large $n$.
Omit at most $r=r(T_o)$ exceptional vertices from $G$, the remaining graph is $G'=G_1'+G_2'$ and each $G_i'=k_iH_i$. Since $n$ is sufficiently large, each block $H_i$ is $K_1$, for otherwise $G_i$ contains $e(T)\cdot K_2$ which is impossible since $e(T)\cdot K_2\in \mathcal{M}(T_o)$. That is, $G'=G_1'+G_2'$ is a complete bipartite graph.  Let $A_i=V(G_i')$ for $i=1,2$.
On the other hand, by the symmetry condition in the definition of $\mathcal{D}(n,2,r)$, any exceptional vertex is adjacent to all vertices in $A_1\cup A_2$, or to all vertices in $A_i$ but no vertices in $A_{3-i}$, or to no vertices in $A_1\cup A_2$.
Let $W$ and $W'$ be the sets of vertices in $V(G)-V(G')$ that are adjacent to all vertices and to no vertex in $G'$, respectively, and $B_i$ be the set of vertices in $V(G)-V(G')$ that are adjacent to all vertices in $A_{3-i}$ but not any vertex in $A_{i}$, $i=1,2$.
If $W'\not=\emptyset$, say $v\in W'$, then we can delete the edges between $v$ and all other exceptional vertices, and then add all edges between $v$ and $A_2$. This does not decrease the number of edges because $r$ is a constant and $A_2\sim \frac{n}{2}$. Also $G$ still contains no $T_o$, for otherwise replace $v$ with some vertex in $A_1$,  we can find a $T_o$ in the original graph $G$. Hence, $W'=\emptyset$.

\begin{Claim}\label{claim1}
 $|W|=a-1$ and $G[W]$ is $\mathcal{B}(T_o)$-free.
\end{Claim}
\pf If $|W|\ge a$, then we can find a copy of $T$ in $G[W\cup A_1]$. By Lemma \ref{lemma1}, $T\in\mathcal{M}(T_o)$ and hence $T_o\subseteq G[W\cup A_1\cup A_2]$, a contradiction. If $|W|\le a-2$, then
\[\begin{split}
 \ e(G)&\le e(G(n,2,a-1))+\binom{|W|}{2}+\binom{|B_1|}{2}+\binom{|B_2|}{2}
\\&=\frac{n^2}{4}+\frac{a-2}{2}n+o(1)<\frac{n^2}{4}+\frac{a-1}{2}n+o(1)=e(G(n,2,a)),
\end{split}\]
a contradiction. Therefore, $|W|=a-1$.

If $G[W]$ is not $\mathcal{B}(T_o)$-free, then there is some $T'$ in $\mathcal{M}(T_o)$ such that $T'\subseteq G[W\cup A_1]$ by  the definition  of $\mathcal{B}(T_o)$,   and hence $T_o\subseteq G[W\cup A_1\cup A_2]$, a contradiction.  $\hfill\square$
\vskip 2mm
Let $G^*=G-W$, $X_i=A_i\cup B_i$ and $G_i^*=G[X_i]$ for $i=1,2$, and $G_{cr}^*=G[X_1,X_2]$.
Moreover,
Let $T=T[A,B]$ be a tree,  $u\in A$ with $d(u)=k$ and $N(u)=\{u_1,u_2,\ldots, u_k\}$. Assume $d_I(u)=k_1$ and $uu_i$ is of Type \Rmnum{1} for $i\le k_1$. Obviously, $k\geq k_1$.  We  now show the graph $G^*$ with partition $V(G^*)=X_1\cup X_2$ satisfying the conditions of Lemma \ref{lemma4} by the following three claims.

\begin{Claim}\label{claim2}
$G_i^*$ is $\{S_{k-\ell}\cup \ell K_2 : 0\leq \ell\leq min\{k-k_1,k-2\}~or~\ell=k-1\}$-free, $i=1,2$.

\end{Claim}
\pf Suppose to the contrary that $G_i^*$ contains an $S_{k-\ell}\cup \ell K_2$ for some $\ell$.

If $\ell=k-1$, then $S_{k-\ell}\cup \ell K_2=kK_2$. Let $T'$ be forest obtained from $T$ by splitting $u$ first and then peeling off the {\it leaf-edges} of Type \Rmnum{2} which are incident with $u$ in $T$. Then $T'$ is the union of $kK_2$ and a forest $T'[A-\{u\},B]$.  Embed $A-\{u\}$ into $W$ and $B$ into the vertices other than that of the $kK_2$ in $A_i\subseteq X_i$, we can find  a $T'$  in $G^*[W\cup X_i]$. By Lemma \ref{lemma1},  $T'\in \mathcal{M}(T_o)$, a contradiction. Hence, $G_i^*$ is $kK_2$-free.

Assume that $0\leq \ell\leq min\{k-k_1,k-2\}$. Let $T'$ be a forest obtained from $T$ by splitting the set $\{u_1,u_2,\ldots, u_k\}$, and then peeling off $\ell$ edges of Type II from the star $S_k$ with center $u$. Then $T'$ is the union of $S_{k-\ell}\cup \ell K_2$ and a forest $T'[A-\{u\},B']$. Embed $A-\{u\}$ into $W$ and $B'$ into the vertices other than that of the $S_{k-\ell}\cup \ell K_2$ in $A_i\subseteq X_i$, we can get a  $T'$  in $G^*[W\cup X_i]$. By Lemma \ref{lemma1},  $T'\in \mathcal{M}(T_o)$, a contradiction. So, $G_i^*$ is $S_{k-\ell}\cup \ell K_2$-free. $\hfill\square$

\begin{Claim}\label{claim3}
For each $i=1,2$ and any vertex $v\in X_i$,
$$d_{X_i}(v)+\nu(G^*[N_{X_{3-i}}(v)])\le k-1.$$
\end{Claim}
\pf By symmetry, we assume  $i=1$. Let $N_{X_1}(v)=\{v_1,\ldots,v_t\}$ and $\{x_1y_1,\ldots,x_\ell y_\ell\}$ be a maximum matching in $G[N_{X_2}(v)]$.
If $t+\ell\ge k$, then note that $G^*[W,A_2]$ is a complete bipartite graph and $|A_2|$ is sufficiently large, we can find a copy of $T$ in $G^*$ by embedding $A-u$ into $W$, $u$ into $v$, $\{u_1,\ldots, u_{k_1}\}$ into $\{v_1,\ldots,v_t\}\cup \{x_1,\ldots, x_\ell\}$ and all other vertices of $B$ into $A_2$.

For each $xy\in E(T)$ with $x\in W$ and $y\in A_2$, since $G^*[A_1,A_2]$ is complete bipartite graph  and
$|A_i|$ is sufficiently large, we can use some vertices in $A_1\cup A_2$ to form an odd cycle $C(x,y)$. For each $vy\in E(T)$,
if $y\in  \{v_1,\ldots,v_t\}$, then
choose a vertex $y'\in A_2$,  and if $y=x_j$ for some $j$ with $1\leq j\leq \ell$, then let $y'=y_j$, we can get a triangle $C(v,y)=vyy'$; If $y\in \{u_{k_1+1},...,u_k\}$, then since both $G^*[A_1,X_2]$ and  $G^*[A_2,X_1]$ are complete bipartite graphs and $|A_i|$ is sufficiently large,  we can get an odd cycle $C(v,y)$ by using one of the $t+\ell-k_1\geq k-k_1$ edges in $\{vv_1,...,vv_t, x_1y_1,...,x_\ell y_\ell\}$, which are not used to form a triangle $vyy'$, together with some vertices in $A_1\cup A_2$.
Thus, $G^*$ contains a $T_o$,  a contradiction. $\hfill\square$

\begin{Claim}\label{claim4}
If $k>k_1$, then $N_{X_{3-i}}(v)$ are isolated in $G^*_{3-i}$ for any $v\in X_i$ with $d_{X_i}(v)=k-1$.
\end{Claim}
\pf By symmetry, we assume $i=1$.  Let $N_{X_1}(v)=\{v_1,\ldots, v_{k-1}\}$  and $zz'\in E(G_2^*)$ with $z\in N(v)$.  Since $k>k_1$ implies $k_1\leq k-1$, $G^*[W,A_2]$ is a complete bipartite graph and $|A_2|$ is sufficiently large, we can find a copy of $T$ in $G^*$ by embedding $A-\{u\}$ into $W$, $u$ into $v$, $\{u_1,...,u_{k_1}\}$ into $\{v_1,...,v_{k_1}\}$ and all other vertices of $B$ into $A_2$.

For each $xy\in E(T)$ with $x\in W$ and $y\in A_2$, since $G^*[A_1,A_2]$ is complete bipartite graph  and
$|A_i|$ is sufficiently large, we can use the vertices in $A_1\cup A_2$ to form an odd cycle $C(x,y)$. For each $vy\in E(T)$, if $y\in  \{v_1,\ldots,v_{k_1}\}$, then choose an unused vertex $y'\in A_2$, we can get a triangle $C(v,y)=xyy'$, and if $y\in A_2$, then since both $G^*[A_1,X_2]$ and  $G^*[A_2,X_1]$ are complete bipartite graphs and $|A_i|$ is sufficiently large, we can get an odd cycle $C(v,y)$ by using one of the $k-k_1$ edges in  $\{vv_{k_1+1},\ldots,vv_{k-1}, zz'\}$ and some other vertices in $A_1\cup A_2$. Hence, $G^*$ contains a $T_o$,  a contradiction. $\hfill\square$

\vskip 2mm
By the structure of $G$, we have
$$e(G)\le e(G[W])+|W||G^*|+e(G^*)
= e(G[W])+|W||G^*|+e(G_1^*)+e(G_2^*)+e(G_{cr}^*).$$
By Claims \ref{claim2}, \ref{claim3} and \ref{claim4}, and applying Lemma \ref{lemma4} on $G^*$, we have
$$e(G_1^*)+e(G_2^*)+e(G_{cr}^*)\le |X_1||X_2|+\begin{cases}
(k-1)^2, & if ~k> k_1,\\
f(k-1,k-1), & if~k=k_1.
\end{cases}$$
By Claim \ref{claim1},  $e(G[W])\le ex(a-1, \mathcal{B}(T_o))$. Moreover, $|W||G^*|+|X_1||X_2|\leq e(G(n,2,a))$ by the definition of the graph $G(n,2,a)$. Therefore, we have
$$e(G)\le e(G(n,2,a))+ex(a-1,\mathcal{B}(T_o))+\begin{cases}
(k-1)^2, &if~ k> k_1,\\
f(k-1,k-1), &if~ k=k_1.
\end{cases}$$
The proof of Theorem \ref{Thm1} is complete.$\hfill\blacksquare$

\section{Remark on Conjecture \ref{conj1}}

It is clear that, if we color the edges of an
extremal graph for $ex(n,H)$ red and color the other edges blue in a $K_n$, then any red edge is not in monochromatic $H$, which  implies $f(n,H)\geq ex(n,H)$ for any $H$. Keevash and Sudakov \cite{Sudakov}  proved that if $H$ has an edge $e$ such that $\chi(H-e)=\chi(H)-1$ or  $H=C_4$, then $f(n,H)=ex(n,H)$. Later, Ma \cite{Ma} and Liu et al. \cite{Liu2} confirmed Conjecture \ref{conj1} for a large family of bipartite graphs, including cycles and some complete bipartite graphs. Recently,  Yuan \cite{Yuan} found some counterexamples to Conjecture \ref{conj1} with large chromatic number. However, it remains unknown  if Conjecture \ref{conj1} holds for all bipartite graphs or other graphs with small  chromatic number. In particular, is it true $f(n,H)=ex(n,H)$ for $\chi(H)=3$?


Let $T_o$ be a good odd-ballooning of $T$. By Theorem \ref{Thm1}, $ex(n,T_o)$ is the sum of three terms, and if $d_I(u)<k$, then $G(n,2,a,\mathcal{B}(T_o), K_{k-1,k-1})$ is an extremal graph. Note that the size of the maximal $\mathcal{B}(T_o)$-free graph embedded into the set $X$ is the middle term $ex(a-1, \mathcal{B}(T_o))$. If we color the edges of such an extremal graph for $ex(n,T_o)$ red and color the other edges blue in a $K_n$, one can see the blue edges in the set $X$ are not covered by any monochromatic $T_o$, too.
Hence, we have

$$f(n, T_o)\ge e(G(n,2,a))+\binom{a-1}{2}+
\begin{cases}
(k-1)^2, & if ~k> k_1,\\
f(k-1,k-1), & if~k=k_1.
\end{cases}$$

Compare the right hand of this inequality with $ex(n,T_o)$ in Theorem \ref{Thm1}, we can see if $ex(a-1,\mathcal{B}(T_o))\not= \binom{a-1}{2}$, then Conjecture \ref{conj1} is not true for $T_o$. By Lemma \ref{lemma2}, we know if $\beta(T)<a$, then $\mathcal{B}(T_o)\not= \{K_a\}$ and hence $ex(a-1,\mathcal{B}(T_o))\not= \binom{a-1}{2}$. It is easy to see there are many trees $T=T[A,B]$ with $\beta(T)<a$, for example, a double star. Let $T=S_{a, b}$ be a double star and $u,v$ two centers of $T$ with $d(u)=a$, $d(v)=b$ and $a\leq b$.
For any good odd-ballooning $T_o$,  because $S_{a,b}\in \mathcal{M}(T_o)$ by Lemma \ref{lemma1} and $\{u,v\}$ is a covering of $S_{a,b}$,  we have $uv\in \mathcal{B}(T_o)$ by the definition of $\mathcal{B}(T_o)$. Therefore,   $ex(a-1,\mathcal{B}(T_o))=0$ and
$$ex(n, T_o)=e(G(n,2,a)).$$
That is to say, any good odd-ballooning of a double star $S_{a,b}$ is a  counterexample to Conjecture \ref{conj1},  and so $f(n,H)>ex(n,H)$ for many $H$ with $\chi(H)=3$.

\vskip 5mm
\noindent{\bf\large Acknowledgements}
\vskip 3mm

This research was supported by NSFC under grant numbers  11871270, 12161141003 and 11931006.


\end{document}